\documentclass[12pt,reqno]{amsart}
\usepackage{amssymb}
\usepackage{mathrsfs}
\usepackage{dsfont}
\usepackage{amsmath,amsthm}
\usepackage{amsfonts}
\usepackage{enumerate, xspace}
\usepackage{a4wide}
\usepackage{graphicx}
\RequirePackage{hyperref}
\usepackage{hyperref}

\theoremstyle{definition}
\newtheorem{thm}{\bf Theorem}[section]

\newtheorem{defin}[thm]{\bf Definition}
\newtheorem{lemma}[thm]{\bf Lemma}
\newtheorem{coro}[thm]{\bf Corollary}
\theoremstyle{remark}

\numberwithin{equation}{section}

\newcommand{\dis}{\displaystyle}
\newcommand{\eps}{\epsilon}
\newcommand{\nn}{\nonumber}
\newcommand{\ga}{\gamma}
\newcommand{\Om}{\Omega}
\newcommand{\om}{\omega}

\newcommand{\und}{\underline}
\newcommand{\la}{\lambda}

\usepackage{color}
\definecolor{titti}{rgb}{0.09, 0.75, 0.50}
\definecolor{darkspringgreen}{rgb}{0.09, 0.45, 0.27}
\definecolor{greyd}{cmyk}{0,0,0,0.4}
	\definecolor{applegreen}{rgb}{0.55, 0.71, 0.0}
	\definecolor{darkpastelgreen}{rgb}{0.01, 0.75, 0.24}
\definecolor{electricpurple}{rgb}{0.75,0.00,1}

\newcommand{\be}[1]{\begin{equation}\label{#1}}
\newcommand{\ee}{\end{equation}}

\let\oldcite\cite
\renewcommand*{\cite}[1]{{\color{blue} \oldcite{#1}}}

\title{Scaling limit of a generalized contact process}
\author{by Logan Chariker$^1$, Anna De Masi$^2$, Joel L. Lebowitz$^3$ and Errico Presutti$^4$,}
\thanks{\newline $^4$Gran Sasso Science Institute, 67100 L'Aquila, Italy. \url{errico.presutti@gmail.com} \newline $^2$Universit\`a degli Studi dell'Aquila, 67100 L'Aquila, Italy. \url{demasi@univaq.it} \newline $^1$School of Natural Sciences, Institute for Advanced Study, Princeton, New Jersey. \url{chariker@ias.edu} \newline $^3$Departments of Mathematics and Physics, Rutgers University, New Brunswick, New Jersey. \url{lebowitz@math.rutgers.edu}
 \newline\newline Keywords: neurons with discrete voltage, integrate and fire, generalized contact process, mean field, spatial dependence, hydrodynamic limit}

\begin{document}

\maketitle

\noindent{\bf Abstract.} We derive macroscopic equations for a generalized contact process that is inspired by a neuronal integrate and fire model on the lattice $\mathbb{Z}^d$.  The states at each lattice site can take values in $0,\ldots,k$.  These can be interpreted as neuronal membrane potential, with the state $k$ corresponding to a firing threshold.  In the terminology of the contact processes, which we shall use in this paper, the state $k$ corresponds to the individual being infectious (all other states are noninfectious).  In order to reach the firing threshold, or to become infectious, the site must progress sequentially from $0$ to $k$.  The rate at which it climbs is determined by other neurons at state $k$, coupled to it through a Kac-type potential, of range $\gamma^{-1}$.  The hydrodynamic equations are obtained in the limit $\gamma\rightarrow 0$.  Extensions of the microscopic model to include excitatory and inhibitory neuron types, as well as other biophysical mechanisms, are also considered.

\section{Introduction}
The derivation of macroscopic deterministic time evolution
equations from underlying microscopic dynamics is one of the central problems of non-equilibrium statistical
mechanics. This micro-to-macro transition is a very difficult mathematical problem with only limited progress so far \cite{Spohn, DP, KL}. This can be overcome to some extent when the underlying microscopic dynamics is stochastic with very strong ergodic properties. Examples are the time evolution of the stochastic Ising model via Glauber or Kawasaki dynamics. There one has derived rigorously macroscopic equations in a space-time scaling limit \cite{GLP,DGLP,GL}. These equations are of the mean field type using long range Kac type interactions on the microscopic scale.

   In this note we derive macroscopic equations presented and partially solved in \cite{CL}. The microscopic model system described here is inspired by neuronal integrate-and-fire models \cite{Gerstner}. In the simple version of this model the membrane voltage increases until it reaches a maximum threshold value at which time it fires (spikes). When it fires that neuron's membrane voltage gets reset to its minimum value. At the same time other neurons connected to it, whose potential is below threshold, increase their potential at a rate depending on the strength of their connectivity to the neuron which has just spiked. In the macroscopic equations we considered in \cite{CL} we discretized the values which the membrane potential can
take restricting it to the integer set $\{0, 1, ......,k\}$. When and only when a neuron is in state $k$, its maximum value, it causes other neurons connected to it with potential $j<k$ to transit to the next level $j+1$. Independently, neurons with potential values $k$ spike and assume the value 0. The neurons in the microscopic model live on the d-dimensional lattice $\mathbb{Z}^d$ with spacing $\gamma$ and their interaction is given by a Kac type function, $J(\gamma |x-y|)$ \cite{GLP}. In the limit $\gamma \rightarrow 0$ one obtains the macroscopic equations.

It turns out that for $k=1$ the model is equivalent to the well known contact process with the state $j=0$,  corresponding to the healthy state and the state $j=1$ the infected one\cite{Liggett,Marro}. For $k>1$ the model can be thought of as a generalized contact process with only the state $j=k$
  being infectious. In terms of neural models the case $k=1$ corresponds to the stochastic Wilson-Cowan model \cite{Goychuk,Zankoc}, which is a popular simplified model of neural systems.  Setting $k>1$ introduces inactive states which behave like subthreshold neuron potentials and leads to more complicated behavior.

  The analysis in this note will be done entirely in the context of the generalized contact process.  We consider an extension of the classical contact process where the state of an individual is described by a potential $U$: when $U=0$ the individual is healthy, when $0<U<k$ is sick but not contagious and when $U=k$ is both sick and contagious (in the classical contact process $k=1$).  Infections are long range and described by a Kac potential with range $\gamma^{-1}$.  We study the system in the macroscopic limit $\gamma \rightarrow 0$.

In more realistic models there are two types of neurons, excitatory ones which act as those described above and inhibitory ones \cite{Kandel,Gerstner}. The latter ones also have a threshold for firing but instead of increasing the potential of other neurons when firing it decreases them. This, as well as other generalizations can be incorporated in the microscopic model studied here. They lead to more complicated macroscopic equations, which we are currently exploring.  Their derivation uses the same formalism as the derivation given here.  These will be discussed briefly at the end of this note.

The outline of the rest of the paper is as follows.  In section 2 we give a precise definition of the microscopic model and present the hydrodynamic limit equations.  In sections 3-7 we prove the hydrodynamic limit.  In section 8 we describe some generalizations of the model.

\section{The model }
\label{sec.n2}

\medskip
{\bf The macroscopic region.}

\noindent
The macroscopic region $\Om$ is a torus in $\mathbb R^d$  of side $L$,  for simplicity  $L$ is a large positive integer.

\medskip
{\bf The microscopic region.}  Let $\ga=2^{-n_1}$, $n_1$ a positive integer.  The microscopic region is the torus $\Om_\ga=\ga^{-1}\Om\cap \mathbb Z^d$.

\noindent

{\bf The Kac potential.}

\noindent
Dynamics is defined in terms of the Kac potential   $J_\ga(x,y) = a_\ga \ga^d J(\ga x,\ga y)$, $x,y\in\Om_\ga$.  $a_\ga$ is the normalization coefficient which makes $J_\ga(x,y)$ a probability; $J(r,r')$ is a smooth, non negative, symmetric probability kernel with finite range $R$.  
 Call
   \begin{equation}
	\label{n2.4bis}
R_{\ga} =\sup\{ |x-y| : J_\ga(x,y)>0\},\quad R_\ga=\ga^{-1}R
	\end{equation}
the range of the interaction $J_\ga(x,y)$.

\bigskip
\noindent The macroscopic limit is defined by letting $\ga\to 0$.

\medskip
{\bf The time evolution.}

\noindent
Time evolution is described by a jump Markov process where there are two types of jumps related to infection and recovery. The individual at
 site $x$ with $U(x)=k$ recovers at rate $1$ and the potential after recovery becomes $U(x)= 0$.
Moreover the individual at
 site $x$ with $U(x)=k$
infects  the one at site $y$ if $U(y)<k$ at rate $\la^* J_\ga (x,y)$
and the effect of the infection is that $U(y) \to U(y)+1$.

We denote by $U_t(x)$, $x\in\Om_\ga$ the potentials at time $t$ and we denote by $P^\ga$ the law of this process in $\Om_\ga$.
%

\medskip

\noindent
\begin{defin}
\label{defin3.0.1}  {\bf The initial condition.}
For any fixed $\ga$ 
the potentials $U_0(x)$, $x\in\Om_\ga$, at time 0 are mutually independent and the following holds:

\medskip
\begin{itemize}

\item   $P^\ga[U_0(x)=i] =\rho_0(\ga x,i)$, $i=0,..,k$ where   $\rho_0(r,i)>0$, $i=0,..,k$ is such that
$\sum_i \rho_0(r,i)=1$ for any $r\in \Om$. The  variables $\rho_0(r,i)$ are the initial densities.

\item The densities $\rho_0(r,i)$, $r\in \Om$, $i=0,..,k$, have values in $[\eps,1-\eps]$, $\eps>0$ and
 are smooth.


\end{itemize}

\end{defin}

Assuming propagation of chaos namely that the potentials $U_t(x)$, $x\in\Om_\ga$ at any $t>0$ are mutually independent and $P^\ga[U_t(x)=i] =\rho_t(\ga x,i)$, $i=0,..,k$, then $\rho_t$ satisfies
 \begin{equation}
	\label{1.2a}
\frac{d}{dt}\rho_t(r,i)=\{ \rho_t(r,i-1) -\rho_t(r,i)\} \int_{\Om}dr'\la^*\rho_t(r',k) J(r',r)
	\end{equation}
when  $1\le i \le k-1$; for $i=k$
\begin{equation}
	\label{1.2b}
\frac{d}{dt}\rho_t(r,k)= - \rho_t(r,k)+ \rho_t(r,k-1)\int_{\Om}dr'\la^*\rho_t(r',k)J(r',r)
	\end{equation}
and  for $i=0$:
\begin{equation}
	\label{1.2c}
\frac{d}{dt}\rho_t(r,0)=  \rho_t(r,k)- \rho_t(r,0)\int_{\Om}dr'\la^*\rho_t(r',k) J(r',r)
	\end{equation}
the initial condition being $\rho_0(r,i)$, $\rho_0$ as in Definition \ref{defin3.0.1}.  These are the macroscopic equations analyzed in \cite{CL}.

 We will use three scaling parameters $\ga=2^{-n_1}$, $ n_1 \in \mathbb N$, $\delta=2^{-n_2}$,  $ n_2 \in \mathbb N$ and $\xi=2^{-n_3}$,  $n_3 \in \mathbb N$. The macroscopic limit, namely the solution of equations \eqref{1.2a}-\eqref{1.2c}, is obtained by first letting $\ga\to 0$, then $\delta\to 0$ and finally $\xi\to 0$.

To prove this limiting behavior we first prove the hydrodynamic limit for a modified dynamics, called the auxiliary process, with a Kac potential $A_\ga$ which is a coarse grained version of $J_\ga$. We then   obtain in the limit $\ga\to 0$ a macroscopic equation with a kernel  $A_\xi$ which is a coarse graining version of $J$, see Theorem \ref{thma1.3bis}. In the limit $\xi\to 0$ we get \eqref{1.2a}-\eqref{1.2c}, see Theorem \ref{thm2.3a}.

\bigskip

\section {The auxiliary process}
\label{sec2.1aa}

\medskip
\noindent

The auxiliary process is defined as the previous one but with a
 piecewise constant  kernel $A_\ga(x,y)$  in the place of $J_\ga(x,y)$.

\medskip

In order to define $A_\ga$ we need the following definition.
%
%
%

\begin{defin} {\em The basic partition}.

\label{defin3.0}
We call
  \begin{equation}
	\label{3.0}
\Om^*_\ga = \{r \in \Om: r=\ga x, x \in \mathbb Z^d\},\quad \Om^*=\bigcup_\ga \Om^*_\ga,\qquad \Om_\ga=\ga^{-1}\Om^*_\ga
	\end{equation}

{\em The basic partition} of $\Om$ is denoted by  $\pi_\xi$, $\xi=L2^{-n_3}$;  its atoms $C^\xi\in\pi_\xi$ are cubes of side $\xi$ and $C^\xi(r)$, $r\in \Om$, denotes the atom  which contains $r$.

The microscopic basic partition  $\pi_{\ga,\xi}$ is made by atoms $C^{\ga,\xi}= \ga^{-1}C^\xi$, $C^\xi\in\pi_\xi$, $C^{\ga,\xi}$ has $N=(\ga^{-1}\xi)^d$ elements.
We write $C^{\ga,\xi}(x) = \ga^{-1}C^\xi(\ga x)$, $C^\xi\in\pi_\xi$.

We say that two  atoms $C^{\ga,\xi}$ and $ D^{\ga,\xi}$ of the basic partition interact with each other if there are $x\in C^{\ga,\xi}$ and $y\in D^{\ga,\xi}$ such that $J_\ga(x,y) >0$.

To simplify the notation we will  drop the superscript $(\ga,\xi)$ from  the cubes $C^{\ga,\xi}$ unless confusion may arise.

\end{defin}

\medskip
{\bf The new piecewise constant kernel.}
In the new process the rate at which $x$ infects $y\ne x$ is
$\la^*A_{  \ga}(x,y)$.
The new kernel $A_\ga(x,y)$  is defined by averaging $J_\ga(x,y)$ over the atoms of the
basic partition $\pi_{\ga,\xi}$, more precisely

  \begin{equation}
	\label{n2.4n}
A_{\ga}(x,y)=
\frac{1}{N^2}\sum_{x'\in  C^{\ga,\xi}(x)}\sum_{y'\in  D^{\ga,\xi}(y)}
 J_\ga(x',y')
	\end{equation}
where
$ J_\ga(x,x)=0$ and  $N=|C^{\ga,\xi}|$.

{\bf Properties of $A_{\ga}$.}

\nopagebreak
$\bullet$\; For any $x$, $A_{\ga}(x,\cdot)$ is a probability, in fact
  \begin{eqnarray}
	\nn
&&\hskip-1cm
\sum_{ y } A_{\ga}(x ,y) = \sum_{D\in \pi_{\ga,\xi}} \sum_{y\in D} A_{\ga} (x,y) =\sum_{D\in \pi_{\ga,\xi}}  A_{\ga} (x,y_D) |D|
=  \frac{1}{N}\sum_{x'\in  C(x)}\sum_{D\in \pi_{\ga,\xi}} \sum_{y'\in  D}
 J_\ga(x',y')
 \\&&\hskip3cm = \sum_{D\in \pi_{\ga,\xi}} \sum_{y'\in  D} J_\ga(x,y')=\sum_{z} J_{\ga} (x,z)=1	
 	\label{n2.5.11}\end{eqnarray}
where $y_D$ is any point in $D$. In the third equality we have used that $|D|=N$.

$\bullet$\; $ A_{\ga}(x ,y)=0$ if $C(x) $ and $ D(y) $ do not interact.

$\bullet$\;
If  instead $C(x) $ and $ D(y) $   interact
  \begin{equation}
	\label{n2.6}
A_{ \ga}(x,y) \le \ga^d   \| J(r,r'|\|_\infty   = :c\xi^d \frac 1N
	\end{equation}

 \medskip

We denote by $P^{A_\ga,\ga}$ the law of the auxiliary process with initial conditions  unchanged, see Definition \ref{defin3.0.1} and by $\mathcal P^{A_\ga,\ga}$ the corresponding  law of
  of the density variables
\begin{equation}
	\label{n2.2}
v^{\ga,\xi}_{i,t} (x)= \frac 1N |C^{\ga,\xi}_{i,t}(x) |,\quad C^{\ga,\xi}_{i,t} = \{ x' \in C^{\ga,\xi}(x): U_t(x')=i\}
	\end{equation}

The following Theorem will be proved  in Section \ref{sec55}.

\begin{thm}
\label{thma1.3bis}
Fix $T>0$ and $t\in[0,T]$, we
denote by $\mathcal P^{A_\ga,\ga}_{t}$ the restriction of $\mathcal P^{A_\ga,\ga}$ to time $t$. In analogy with \eqref{n2.4n} we define the kernel $A_\xi(r,r')$, $r,r' \in \Om^*$ as
  \begin{equation}
	\label{nuovo}
 A_\xi(r,r') = \frac{1}{\xi^{2d}}\int_{  C^{\xi}(r)}dr_1\int_{ C^{\xi}(r')}dr'_1
 J(r_1,r_1')
  \end{equation}
 Then
$\mathcal P^{A_\ga,\ga}_{t}$  converges  as $\ga \to 0$ to a probability $\mathcal P^{A_\xi}_{t}$ which is supported by $\varphi_\xi(r,i;t)$ $r\in \Om$ which is the solution
at time $t$ of  the equations
 \begin{equation}
	\label{n2.3.3.0}
\frac{d}{dt}\varphi_\xi(r,i;t)=\{\varphi_\xi(r,i-1;t) -\varphi_\xi(r,i;t)\} \int_{\Om}dr'\la^*\varphi_\xi(r',k;t) A_\xi(r',r)
	\end{equation}
when  $1\le i \le k-1$; for $i=k$
\begin{equation}
	\label{n2.3.3.0bis}
\frac{d}{dt}\varphi_\xi(r,k;t)= -\varphi_\xi(r,k;t)+ \varphi_\xi(r,k-1;t)\int_{\Om}dr'\la^*\varphi_\xi(r',k;t)A_\xi(r',r)
	\end{equation}
and  for $i=0$:
\begin{equation}
	\label{n2.3.3.0ter}
\frac{d}{dt}\varphi_\xi(r,0;t)=  \varphi_\xi(r,k;t)- \varphi_\xi(r,0;t)\int_{\Om}dr'\la^*\varphi_\xi(r',k;t) A_\xi(r',r)
	\end{equation}
the initial condition being $\dis{\varphi_\xi(r,i;0)=\frac 1{|C^\xi(r)|} \int_{C^\xi(r)}dr'\;\rho_0(r',i)}$, $\rho_0$ as in Definition \ref{defin3.0.1}.

\end{thm}

\medskip
\noindent
Observe that $\varphi_\xi(r,i;t)$, $r\in \Om$  is constant on the cubes $C^\xi$.
Convergence of the space-time joint distribution of the densities will be proved in Section \ref{sec5} together with the following Theorem.

\medskip
\begin{thm}
\label{thm2.3a}
The function $\varphi_\xi(r,i;t)$, $r\in \Om$ , $t\in[0,T]$ converges as $\xi\to 0$ to $\rho_t(r,i)$ solution of the equations \eqref{1.2a}-\eqref{1.2c}.
\end{thm}

\medskip
\centerline{Sketch of the proof of Theorem \ref{thma1.3bis}}

\nopagebreak

\noindent
In the theory of hydrodynamic limit for stochastic interacting particle systems
a typical procedure is to use the martingale decomposition for the variables of interest, see for instance the book  \cite{KL}.  Applied to our case we have
\begin{equation}
	\label{sketch.1}
  v^{\ga,\xi}_{i,t} (x) -  v^{\ga,\xi}_{i,0} (x)
  = \int_0^t ds\, L_\ga v^{\ga,\xi}_{i,s} (x) + M^{\ga,\xi}_{i,t} (x)
	\end{equation}
where $ L_\ga$ is the generator of the process and $M^{\ga,\xi}_{i,t} (x)$ is a martingale. $M^{\ga,\xi}_{i,t} (x)$ is a ``fluctuation term'' and one can often prove that in the hydrodynamic limit $N\to \infty$ $M^{\ga,\xi}_{i,t} (x)$ vanishes with probability going to 1.  The hardest problem is to control $L_\ga v^{\ga,\xi}_{i,t} (x)$ whose explicit expression for $1\le i \le k-1$ in our case is
\begin{equation}
	\label{sketch.2}
L_\ga v^{\ga,\xi}_{i,t} (x)= \Big(v^{\ga,\xi}_{i-1,t} (x)-v^{\ga,\xi}_{i,t} (x)\Big) \la^* \sum_y A_\ga(y,x) v^{\ga,\xi}_{k,t}(y)
	\end{equation}
By compactness $v^{\ga,\xi}_{i,t} (x)$ converges (by subsequences)  weakly in probability to some limit density but the problem is that in \eqref{sketch.2} the functions $v^{\ga,\xi}$ appear quadratically and in general the weak limit of a product is not the product of the weak limits of the factors.

To close the equations one then needs to prove a factorization property for the $v^{\ga,\xi}_{i,t}$  , i.e.\ propagation of chaos or local equilibrium.  We overcome this difficulty by using the same method as in  \cite{DGLP} and \cite{DOR}.
We discretize time, see  Section \ref{sec.n2bis}:  we use a mesh $\delta$ which will vanish after taking the limit $\ga \to 0$ and study the process in the generic time interval $[n\delta,(n+1)\delta]$ with $n \le \delta^{-1}T$ having conditioned on the values of the potential $U_t(x)$ at time $t=n\delta$.

The increments of the densities in a time interval $[n\delta,(n+1)\delta]$ (having fixed the potentials at time $n\delta$) are given in \eqref{n2.7.1}--\eqref{n2.7.3} in terms of variables $M_{C,D;i}$, $C$, $D$ cubes of the basic partition, $i\in \{0,..,k\}$, and variables $M_{D}$, see \eqref{n2.6.1}.

Probability estimates on  $\mathcal M_{C,D;i}$ are obtained in Theorem \ref{thm.n2.1}, those for $M_{D}$ in Theorem \ref{thm.n2.1.33}. The proof of Theorem \ref{thma1.3bis} is then given in Section  \ref{sec55}.

 The crucial point is  to prove the probability estimates stated in Theorem \ref{thm.n2.1} and in Theorem \ref{thm.n2.1.33}.  We use a graphical representation of the process
where we represent by an arrow $(x,y)$ the infection to the individual at $y$ due to the individual at $x$; the recovery of an individual at $x$ is described  by a ``marked point''.

The collection of arrows and marked points define a natural graph structure, see the paragraph {\it A graph structure} in the next section.
To reconstruct the true process we introduce time variables $\und t(x,y)$ and $\und t(x)$, $\und t(x,y)$ is a finite sequence of times $t_m(x,y)$ and $\und t(x)$ of times $t_m(x)$. $\und t(x,y)$ and $\und t(x)$ are mutually independent Poisson processes with mean $\la^* A_\ga(x,y)$ and respectively mean 1.

The above graph structure is realized by drawing an arrow $(x,y)$ at a time $t\in \und t(x,y)$ and a marked point at $x$ at $t\in\und t(x)$.  Knowledge of all $\und t(x,y)$ and $\und t(x)$ allows to reconstruct the true process, see the paragraph {\em A realization of the process: the clock process} in Section \ref{sec.n2bis}.
 However to know whether at $t=t(x,y)$ there is an infection we need to know all the values $t(x',y';s)$ and $t(x';s)$ for all $s\le t$ as well as the values of the initial potentials.

The analysis of the graph structure of arrows and marked points ignoring
the times when they are drawn is quite simple  because the variables
$\und t(x,y)$ and $\und t(x)$ are mutually independent.  The first crucial point is that an arrow $(x,y)$ corresponds to an infection if at the initial time $U(x)=k$ and $U(y)=i$, $i<k$, provided that the cluster containing $(x,y)$ is made only by the arrow $(x,y)$, see Lemma \ref {lemma.n2.2} and the paragraph  {\em A graph structure} in the next section for the definition of clusters.  Analogous property holds for marked points.  Thus when clusters have only one element the time when the event occurs is not relevant.

The second crucial point is that clusters with more than one element are probabilistically negligible.  An estimate is proved in Corollary
\ref {lemma.n2.2.0}. As argued after \eqref{3.19} this is good enough for clusters with at least 3 elements, for clusters with only two elements we have a more refined argument proved in Lemma \ref {lemt?}.

The crucial step in the proof of  Corollary
\ref {lemma.n2.2.0} is to reduce to a branching process which is studied in  Appendix \ref  {app.A}.


\bigskip

\section{Time discretization and a realization of the process}
\label{sec.n2bis}

{\bf Time discretization.}

\noindent
We discretize time with mesh $\delta =2^{-n_2}$, $n_2\ge 1$.  We fix $\delta$ and a time interval $[n\delta, (n+1)\delta]$,  for a while we will study the process in such a time interval having conditioned on the values $U_{n\delta}$ of the potentials at time $n\delta$.
By choosing $\delta$ small enough the process becomes considerably simpler and we will exploit the following realization of the process.

\medskip

{\bf A realization of the process: the clock process.}

\noindent
We attach to any ordered pair $(x,y)$, $x\ne y$,  independent clocks called $(x,y)$-clocks   which ring
at exponential rate $\la^*A_{ \ga}(x,y)$.  The clocks start  at time $n \delta$ and are stopped at time $(n+1)\delta$, recall that we are studying the process restricted
to the time interval $[n \delta, (n+1) \delta]$. We denote by $\und t(x,y)$ the times when the  $(x,y)$-clock rings.
We introduce also  $x$-clocks  which ring at rate 1, $\und t(x)$ being the times when the $x$-clock rings.  All the above clocks are independent of each other.

The true process is recovered as follows.  If the $(x,y)$-clock rings and at the time of the ring $U(y)<k$ and $U(x)=k$ then $U(y) \to U(y)+1$. Moreover if the $x$-clock rings  at a time when $U(x)=k$ then $U(x)\to 0$. All these rings are {\em effective} while the other rings where the above conditions are not fulfilled   are {\em ineffective}, the potentials are unchanged and they can be ignored.  However it is a very complicated task to understand whether a clock  ring is or is not effective, it depends on all the clock rings $\{\und t(x,y);\und t(x)\}$.  As already mentioned it is convenient to introduce a graph structure.

\medskip

 {\bf A graph structure.}
When the $(x,y)$-clock rings we draw an oriented arrow $(x,y)$, when the $x$-clock rings we draw  a marked point at $x$.
Two arrows are connected if they have a point in common, a marked point is connected to an arrow if it is one of the two points of the arrow.  Clusters are the maximal connected sets of marked points and arrows.  Notice that a same arrow may appear several times in a cluster as well as a same marked point.  We denote by $\mathbf C_1$ the clusters made by a single element, i.e.\ either a marked point or an arrow.
$\mathbf C_j$ are the clusters with $j$ elements.  We will see that  if the time mesh $\delta$ is small the relevant clusters  are the single clusters $\mathbf C_1$. In such a case we have:
\begin{lemma}
\label{lemma.n2.2}
Let $U_{n\delta}(x)=k$, $U_{n\delta}(y)=i$, $i<k$, and let $ \mathbf C_1=(x,y)$ then
$U_{(n+1)\delta}(x)=k$ and $U_{(n+1)\delta}(y)=i+1$.  Analogously if $x$ is a {\em marked point} with $U_{n\delta}(x)=k$ and $ \mathbf C_1=x$ then $U_{(n+1)\delta}(x)=0$
%

\end{lemma}

 \medskip
 \noindent
{\bf Proof.}  The potentials $U_{n\delta}(x)$ and $U_{n\delta}(y)$ can only change when the $(x,y)$-clock rings because $(x,y)\in \mathbf C_1$ and all the other arrows are not connected to $(x,y)$ nor the marked points. Then by the assumption $U_{n\delta}(x)=k$ and $U_{n\delta}(y)=i$ the $(x,y)$-ring is effective hence the statement on the lemma.  The case of $C_1=x$ is proved similarly.   \qed

\medskip

\begin{defin}
\label{def3.5.kkkk}
 In the sequel we will denote by $c$  constants which do not depend on $N$, $\xi$  and $\delta$.

\end{defin}

%
%
%

\begin{thm}
\label{thm?}
For any $a \in (0,1)$ and any $\eps>0$ such that $1-a-2\eps>0$, there is a constant $c$ so that for any $\ga$ and $\delta$ small enough the following holds.  Let $(x,y)$ be an arrow then for any two atoms $ C$ and $ D$  of the basic partition
 \begin{equation}
	\label{n2.6.1.1}
 N^{-1}\sum_{x\in  C}\sum_{y\in  D}\sum_{j\ge 1}\delta^{-aj} \sum_{\mathbf C_j \ni (x,y) }P^{A\ga,\ga}\big[\mathbf C_j\big] \le c\delta^{1-a-2\eps}
	\end{equation}
Moreover 
\begin{equation}
	\label{n2.6.1.2}
 N^{-1}\sum_{x\in C}\sum_{j\ge 1}\delta^{-aj}\sum_{\mathbf C_j\ni x }
 P^{A\ga,\ga}\big[\mathbf C_j \big]
   \le c\delta^{1-a-2\eps} 
	\end{equation}
\end{thm}

\medskip
 We will use the following consequence of
Theorem \ref{thm?}:

\medskip

\begin{coro}
\label{lemma.n2.2.0}
 Let $ C$ and $D$ be two atoms of the basic partition, then for any $j^*\ge 1$
 \begin{equation}
	\label{3.8aa}
P^{A\ga,\ga} [ N^{-1}\sum_{x\in C}\sum_{y\in D}\sum_{j\ge j^*}\mathbf 1_{\mathbf C_j\ni (x,y)}>\delta^{bj^*}]\le c\delta^{(a-b)j^*+1-a-2\eps},\quad a\in(0,1), \quad b\in(0,1)
 	\end{equation}
for  $\eps$ as in Theorem \ref{thm?}.
Analogously
\begin{equation}
	\label{3.8aabb}
P^{A\ga,\ga} [ N^{-1}\sum_{x\in C}\sum_{j\ge j^*}\mathbf 1_{\mathbf C_j\ni x}>\delta^{bj^*}]\le c\delta^{(a-b)j^*+1-a-2\eps},\quad a\in(0,1), \quad b\in(0,1)
 	\end{equation}
\end{coro}

 \medskip
 \noindent
{\bf Proof.}
 We take $a>0$   and  write
\begin{equation*}
 \sum_{j\ge j^*}\sum_{\mathbf C_j \ni (x,y)}P^{A\ga,\ga} [ \mathbf C_j ] = \sum_{j\ge j^*}\delta^{ ja}\delta^{-ja}\sum_{\mathbf C_j \ni (x,y)}P^{A\ga,\ga} [ \mathbf C_j ]
 \le \delta^{ j^*a}\sum_{j}\delta^{-ja}\sum_{\mathbf C_j \ni (x,y)}P^{A\ga,\ga} [ \mathbf C_j ]
	\end{equation*}
Then using Theorem \ref{thm?}  we have
\begin{equation}
	\label{n2.6.1.3}
N^{-1}\sum_{x\in  C}\sum_{y\in  D} \sum_{j\ge j^*}\sum_{\mathbf C_j \ni (x,y)}P^{A\ga,\ga} [ \mathbf C_j ] \le c \delta^{j^*a}\delta^{1-a}
	\end{equation}
By the Markov inequality we then have \eqref{3.8aa}.  The proof of \eqref{3.8aabb} is similar and omitted.\qed

\medskip
For simplicity, in the sequel we write $P^\ga$ instead of $P^{A\ga,\ga}$.
\medskip

%

\begin{defin}
\label{defin3.9.00}
We denote by $\kappa_{x,y,i}(n)$, $i<k$, the number of effective $(x,y)$-rings in the time interval $[n\delta,(n+1)\delta)$, namely those such that when the clock rings $U(x)=k$ and $U(y)=i$; we  denote by $\kappa_{x}(n)$ the number of effective $x$-rings, namely the times $t$ in $\und t(x)$ when $U_t(x)=k$.
We then define for two cubes $C$ and $D$
  \begin{equation}
	\label{n2.6.1}
\mathcal M_{ C, D;i} (n)= \sum_{x\in  C}\sum
_{y\in  D}\kappa_{x,y,i}(n),\,i<k; \quad
\mathcal M_{D}(n)=  \sum_{x\in D}\kappa_{x}(n)
 \end{equation}

\end{defin}
\noindent
Since in the following $n$ is fixed we drop the dependence on $n$ in \eqref{n2.6.1}.

 \vskip.5cm

Recalling \eqref{n2.2} for notation we have for $0<i<k$
   \begin{equation}
	\label{n2.7.1}
v_{i,(n+1)\delta}(y) -v_{i,n\delta}(y)= \frac 1N
\sum_{C}[\mathcal M_{C,D(y);i-1}-\mathcal M_{C,D(y);i}]
	\end{equation}
   \begin{equation}
	\label{n2.7.2}
v_{0,(n+1)\delta}(y)-v_{0,n\delta}(y)
=
\frac 1N \Big(\mathcal M_{D(y)}-
\sum_{C}\mathcal M_{C,D(y);0}\Big)
	\end{equation}
  \begin{equation}
	\label{n2.7.3}
v_{k,(n+1)\delta}(y) -v_{k,n\delta}(y)
=\frac 1N \Big( -
\mathcal M_{D(y)}+
 \sum_{C}
\mathcal M_{C,D(y);k-1}\Big)
	\end{equation}

\medskip

\section{Probabilty estimates}
\label{sec.n2ter}
The aim  is now to get estimates on $\mathcal M_{D}$ and $\mathcal M_{C,D;i}$, $i<k$, see Definition \ref{defin3.9.00}.   They are in general very complicated functions in the space of the clock  rings $\{\und t(x,y);\und t(x)\}$, we shall see however that only cases with few rings are important, the others give a small contribution.  This will be the crucial point in the proof of the following theorem, Theorem
\ref{thm.n2.1}, which  concerns
 the number of events where an individual in a cube $C$ with potential $k$ infects an individual in the cube $D$ which has potential $i<k$.

\medskip

\centerline{\em Probability estimates on  $\mathcal M_{C,D;i}$}

Recall that  we have fixed a time interval $[n\delta,(n+1)\delta]$ and we will not made explicit the dependence on such interval unless confusion may arise.
\begin{thm}
\label{thm.n2.1}


     There are $\theta >0$,  $a\in (\frac 12,1)$,  $b\in(\frac 13,\frac 23a)$ and a constant $c$ so that for all $i<k$ and all $x$ such that $U_{n\delta}(x)=k$ and all $y$ so that $U_{n\delta}(y)=i$
  	\begin{eqnarray}
\nn&&P^\ga\Big[\big|N^{-1}\mathcal M_{ C(x) , D(y); i} -   \la^* \delta N A_\ga(x,y)  v_{k,n\delta}(x) v_{i,n\delta}(y)
 \big|\le cN^{-\eps}+c( \delta N^{-1} + N^{-\theta}+\delta^{3b})\Big]
 \\&&\hskip1cm\ge 1 -c\Big(N^{-\eps} +N^{-1+2\theta}\xi^d\la^*
 +\delta^{1+2a-3b}+c\frac{\delta^{a}}N\Big)
 	\label{n2.2.4}
 \end{eqnarray}
Notice that the conditions $a\in (\frac 12,1)$ and $b\in(\frac 13,\frac 23 a)$ imply that $2a-3b>0$.
\end{thm}

\noindent
The proof will be obtained in several steps.

\medskip
The first step is to reduce to cases where $\und t(x,y)|=1$.  To this end,
recalling Definition \ref{defin3.9.00}, we
consider two cubes  $ C$ and $D$ and for $i<k$ we write
	\begin{equation}
	\label{3.14}
\mathcal M_{ C, D;i} = \sum_{x\in  C}\sum
_{y\in  D}\kappa_{x,y,i}\mathbf 1_{|\und t(x,y|=1}
+R_{ C, D;i}
\end{equation}
where
$$R_{ C, D;i}=
\sum_{x\in  C}\sum_{y\in  D}\kappa_{x,y,i} \mathbf 1_{|\und t(x,y)|>1} $$


\medskip

\begin{lemma}
\label{lemr}
There is a constant $c$ (independent of $N$ and $\delta$) so that
for any $\eps>0$: 
\begin{equation}
	\label{3.16}
P^\ga[\frac 1N R_{ C, D;i}>N^{-\eps}]\le c\frac {\delta^2}{N^{1-\eps}}
\end{equation}

\end{lemma}

\noindent
{\bf Proof. }
Recall that by definition
    \begin{equation}
	\label{3.15}
	P^\ga[|\und t(x,y)|=n] = e^{-\la^*\delta A_\ga (x, y)}\;\frac{\big[\la^*\delta A_\ga (x, y)]^n} {n!}
    \end{equation}
We then bound

\begin{equation*}
	\label{3.17}
\frac 1N R_{ C, D;i}\le \frac 1N \sum_{x\in  C}\sum_{y\in  D}|t(x,y)|\mathbf 1_{|\und t(x,y|>1}
\end{equation*}
and using \eqref{n2.6} we get
\begin{equation*}
E^\ga[\frac 1N R_{ C, D;i}]=\frac 1N \sum_{x\in  C}\sum_{y\in  D}\sum_{k=2}^\infty k e^{-\la^*\delta A_\ga (x, y)}\frac{[\la^*\delta A_\ga (x, y)]^k}{k!}\le c\frac {\delta^2}N
\end{equation*}

By the Markov inequality we then get \eqref{3.16}.\qed

\bigskip

Let $x$ be such that $U_{n\delta}(x)=k$ and $U_{n\delta}(y)=i$, then
	\begin{eqnarray}
	\label{3.18a}
	&&\hskip-1cm\kappa_{x,y,i}\mathbf 1_{|\und t(x,y)|=1}=\mathbf 1_{|\und t(x,y)|=1} \mathbf 1_{\{(x,y)= \mathbf C_1\}}+\kappa_{x,y,i}\mathbf 1_{|\und t(x,y)|=1}
 \sum_{j\ge 2}\sum_{\mathbf C_j}\mathbf 1_{\{(x,y)\in \mathbf C_j\}}
 \\&&\hskip2cm=\mathbf 1_{|\und t(x,y)|=1} +\mathbf 1_{|\und t(x,y)|=1}[\kappa_{x,y,i}-1 ]
 \sum_{j\ge 2}\sum_{\mathbf C_j}\mathbf 1_{\{(x,y)\in \mathbf C_j\}}
 \nn
	\end{eqnarray}
Thus from \eqref{3.14} we get
\begin{equation}
	\label{3.18}
|\frac 1N\mathcal M_{ C, D;i} -\frac 1N \sum_{x\in  C}\sum
_{y\in  D}\mathbf 1_{|\und t(x,y)|=1}|\le \frac 1N R_{ C, D;i}+ \frac 1NT_{ C,  D}
	\end{equation}
where
\begin{equation}
\label{3.19a}
T_{ C, D}=\sum_{j\ge 2} T_{ C, D}^{(j)},\quad T_{ C, D}^{(j)}=
\sum_{x\in  C}\sum_{y\in  D}\mathbf 1_{\{(x,y)\in \mathbf C_j\}}
	\end{equation}
having used that $0\le \kappa_{x,y,i} \le 1$ if
$\mathbf 1_{|\und t(x,y)|=1}$.

\medskip
By  Corollary \ref{lemma.n2.2.0}  we have:

\begin{equation}
	\label{3.19}
P^\ga[\frac 1N  \sum_{j'\ge j}T_{ C, D}^{(j')}>\delta^{jb}]\le c\delta^{(a-b)j+1-a}
\quad a\in(0,1), \quad b\in(0,1)
\end{equation}
We will eventually need to iterate the estimate over all the time intervals $[\delta n,\delta (n+1)]$, i.e.\ $\delta^{-1}$ times,
so that we want $\delta^{-1}\delta^{jb}$ and $\delta^{-1}\delta^{j(a-b)+1-a}$ to vanish when $\delta\to 0$.
When applied to the case $j=2$, the above requires that $b> 1/2$ and also that $b< a/2<1/2$, so that for $j=2$ the conditions cannot be fulfilled.

Instead  if  $j\ge 3$ by
 choosing
$$a\in (\frac 12,1)\quad  \text{and } b\in(\frac 13,\frac 23 a)\quad  \text{ so that }\quad 2a-3b>0$$
we have \begin{equation}
	\label{3.21abc}
-1+jb >0,\quad -1 +j(a-b)+1-a >0,\qquad \forall j\ge 3
	\end{equation}
	\medskip
	
The analysis of $ T^{(2)}_{ C, D}$ requires a more refined estimate which is the content of the next lemma.

\begin{lemma}
\label{lemt?}
\begin{equation}
	\label{3.19.bis}
P^\ga[\frac 1N  T^{(2)}_{ C, D}>\delta^{2-a}]\le c\frac{\delta^{a}}N
\end{equation}

\end{lemma}

\medskip
\noindent
{\bf Proof.}
A cluster  $\mathbf C_2$ which contains the arrow $(x,y)$ is  equal to
   $$
\{(x,y),(y,z)\}\cup \{(x,y),(z,y)\} \cup \{(x,y),(x,z)\}\cup \{(x,y),(z,x)\}
\cup \{(x,y),x)\}\cup \{(x,y),y)\}
   $$
in the last two terms besides the arrow $(x,y)$ the single poins $x$ and $y$ are marked points and in the previous terms $z$ is any point different from $x$ and $y$..

Since the estimates are similar for simplicity we just examine the case with two arrows, $(x,y),(y,z)$.
We  denote by  $\eta_{x,y,z}\in\{0,1\}$ the indicator of this set, thus
     $$
\eta_{x,y,z}=\mathbf 1_{|\und t(x,y)|=1}\mathbf 1_{|\und t(y,z)|=1}
\mathbf 1_{\mathbf C_2 =\{(x,y),(y,z)\}}
     $$
We call
	\begin{equation*}
F=\frac 1N\sum_{x\in  C}\sum
_{y\in  D}\sum_{z}\eta_{x,y,z}
	\end{equation*}
We first compute the expectation recalling that the clocks are independent (see the paragraph {\em A realization of the process: the clock process}) and using \eqref{n2.6} we get
	\begin{equation*}
\mathbb E\big[ F\big]= \frac 1{N^3}\sum_{x\in C}\sum_{y\in D} \sum_{z}
e^{-\la^*\delta (A_\ga( x , y) +A_\ga( y , z))}(\la^*\delta)^2NA_\ga( x , y)
NA_\ga( y , z)) \le c \xi^d\delta^{2}
	\end{equation*}
We next compute the variance and using independence we get
	\begin{eqnarray*}
&&  \hskip-.6cm
E^\ga\Big[(F- E^\ga[F])^2\Big]= \frac 1{N^4}\sum_{x\in C, y\in D} \sum_{z}
e^{-\la^*\delta (A_\ga( x , y) +A_\ga( y , z))}(\la^*\delta)^2NA_\ga( x , y) NA_\ga( y , z)) -E^\ga[F]^2
 \\&& \hskip3cm
 \le \frac 1N c \xi^d\delta^{2}
	\end{eqnarray*}
Thus \eqref{3.19.bis} follows from Chebishev inequality.
\qed

\bigskip

\noindent
{\bf Proof of Theorem \ref{thm.n2.1}.}  We fix $\bar x$ and $\bar y$ as in the hypothesis of the Theorem and we call $ \bar C=C(\bar x)$ and $\bar D=D(\bar y)$.

We write (recall \eqref{3.18})
		\begin{eqnarray}
		\nn
&&\hskip-2cm \Big|N^{-1}\mathcal M_{ \bar C , \bar D; i} -   \la^* \delta N A_\ga(\bar x,\bar y)  v_{k,n\delta}(\bar x) v_{i,n\delta}(\bar y)\Big|\le \big|S_N- \la^* \delta N A_\ga(\bar x,\bar y)  v_{k,n\delta}(\bar x) v_{i,n\delta}(\bar y)\big|
	\\&&\hskip6cm + \frac 1NR_{ \bar C, \bar D;i}+ \frac 1N T_{\bar C,\bar D}
	\label{3.22a}
		\end{eqnarray}
where
\begin{equation}
	\label{3.22b}
 S_N=\frac 1N \sum_{x\in  \bar C}\sum_{y\in \bar D}\mathbf 1_{|\und t(x,y)|=1},
		\end{equation}
The term $R_{ \bar C, \bar D;i}$ has been treated in Lemma \ref{lemr} and the one with $T_{\bar C,\bar D}$ is estimated in \eqref{3.19} for $j=3$ and Lemma \ref{lemt?} for $j=2$.

We  first compute $E^\ga[S_N]$, from \eqref{3.15}
we get
\begin{eqnarray}
	\nn
&&\hskip-1.3cm  E^\ga[S_N]=e^{-\la^*\delta A_\ga (\bar x, \bar y)}\la^*\delta NA_\ga (\bar x, \bar D)   \frac{|\bar C|}{N} \,
\frac{|\bar D|}{N}\\&&= e^{-\la^*\delta A_\ga (\bar x, \bar y)}\la^*\delta N A_\ga (\bar x, \bar D) v_{k,n\delta}(\bar x) v_{i,n\delta}(\bar y)
	\label{3.24bb}
\end{eqnarray}
By \eqref{n2.6} $\dis{NA_\ga (\bar x, \bar D)\le c\la^*\xi^d }$ so that the right hand side of \eqref{3.24bb} is bounded by $\delta \la^*c\xi^d$. Since the clocks are independent we get
	\begin{equation}
	\label{4.14a}
 \mathbb E^\ga\Big[\big(S_N-\mathbb E^\ga[S_N]\big)^2\Big]=\frac 1{N^2} \sum_{x\in  \bar C}\sum_{y\in \bar D}\big\{P^\ga(|\und t(x,y)|=1)-P^\ga(|\und t(x,y)|=1)^2\}\le c\frac {\delta\la^*\xi^d}N
		\end{equation}
By  \eqref{n2.6}
	\begin{equation}
	\label{4.14b}
	1-e^{-\la^*\delta A_\ga (\bar x, \bar y)}\le c \frac \delta N
	\end{equation}
Thus, given any $\theta>0$ by Chebishev inequality  and using\eqref{4.14a}  and \eqref{4.14b} we get
	\begin{eqnarray}\label{3.26}
P^\ga\Big[\big|S_N- \la^* \delta N A_\ga(\bar x,\bar y)  v_{k,n\delta}(\bar x) v_{i,n\delta}(\bar y)\big|>N^{-\theta}\Big]\le c \frac \delta N+c\frac {\delta\la^*\xi^d}{N^{1-2\theta}} 	
     \end{eqnarray}
We thus get the theorem.\qed
\bigskip

\centerline {\it Probability estimates of $\mathcal M_D$}

\medskip
The analysis of $\mathcal M_D$  defined in \eqref{n2.6.1} is very similar to the one we did for $\mathcal M_{C,D;i}$ and sketched below.  The analogue of Theorem \ref{thm.n2.1} is:

\medskip

\begin{thm}
\label{thm.n2.1.33}
 There are $\theta >0$,  $a\in (\frac 12,1)$,  $b\in(\frac 13,\frac 23a)$ and a constant $c$ so that for all  $x$ such that $U_{n\delta}(x)=k$
 	\begin{eqnarray}
\nn&&P^\ga\Big[\big|N^{-1}\mathcal M_{D(x)} -   \delta  v_{k,n\delta}(x)
 \big|\le cN^{-\eps}+c( \delta N^{-1} + N^{-\theta}+\delta^{3b})\Big]
 \\&&\hskip1cm\ge 1 -c\Big(N^{-\eps} +\delta N^{-1+2\theta}
 +\delta^{1+2a-3b}+c\frac{\delta^{a}}N\Big)
 	\label{n2.2.4.43}
 \end{eqnarray}

\end{thm}

\medskip

We write
	\begin{equation}
	\label{3.14.33}
\mathcal M_{D} = \sum_{x\in D}\kappa_{x}\mathbf 1_{|\und t(x|=1}
+R_{D},\quad R_{ D}=
\sum_{x\in D} \kappa_{x} \mathbf 1_{|\und t(x)|>1}
\end{equation}
%
%
Proceeding as in the proof of Lemma \ref{lemr} we have (proof is omitted):

\begin{lemma}
\label{lemr.33} For all $\eps>0$ we have
\begin{equation}
	\label{3.16.35}
P^\ga[\frac 1N R_{ D}>N^{-\eps}]\le c\frac {\delta^2}{N^{1-\eps}}
\end{equation}

\end{lemma}

\medskip
Analagously to  \eqref{3.18a} we write
	\begin{eqnarray}
	\label{3.18a.36}
	&&\hskip-1cm\kappa_{x}\mathbf 1_{|\und t(x)|=1}=\mathbf 1_{|\und t(x)|=1} \mathbf 1_{\{x=\mathbf C_1\}}+\kappa_{x}\mathbf 1_{|\und t(x)|=1}
 \sum_{j\ge 2}\sum_{\mathbf C_j}\mathbf 1_{\{x\in \mathbf C_j\}}
 \\&&\hskip2cm=\mathbf 1_{|\und t(x)|=1} +\mathbf 1_{|\und t(x)|=1}[\kappa_{x}-1 ]
 \sum_{j\ge 2}\sum_{\mathbf C_j}\mathbf 1_{\{x\in \mathbf C_j\}}
 \nn
	\end{eqnarray}
Thus 
\begin{equation}
	\label{3.18.37}
|\frac 1N\mathcal M_{ D} -\frac 1N \sum
_{x\in  D}\mathbf 1_{|\und t(x)|=1}|\le \frac 1N R_{D}+ \frac 1N T_{ D}
	\end{equation}
where
	\begin{equation}
\label{3.19a??.39}
T_{   D} 
=T_{ C, D}^{(2)} +
\sum_{x\in  D} \sum_{j\ge 3}\sum_{\mathbf C_j}\mathbf 1_{\{x\in \mathbf C_j\}}
	\end{equation}

\medskip
The following lemma  is a consequence of Corollary \ref{lemma.n2.2.0} and an argument very similar to the one of Lemma \ref{lemt?}:

\begin{lemma}
\label{lemt.333}
For any  $a\in (\frac 12,1)$ and $b\in(\frac 13,\frac 23 a)$ we have
\begin{equation}
	\label{3.19.40}
P^\ga[\frac 1N  \sum_{x\in  D} \sum_{j\ge 3}\sum_{\mathbf C_j}\mathbf 1_{\{x\in \mathbf C_j\}}>\delta^{3b}]\le c\delta^{1+2a-3b},\quad
P^\ga[\frac 1N  T^{(2)}_{ D}>\delta^{2-a}]\le c\frac{\delta^{a}}N
\end{equation}

\end{lemma}

\noindent{\bf Proof.}
A cluster  $\mathbf C_2$ which contains the marked point $x$ is equal to
   $$
\{x,(x,y)\}\cup   \{x,(y,x) \}
   $$
We only examine the case  $\{x,(x,y)\}$ the other being similar.
%
We call
	\begin{equation*}
F=\frac 1N\sum_{x\in D}\sum_{y\in  D} \eta_{x,y},\qquad  \eta_{x,y}=\mathbf 1_{|\und t(x,y)|=1}\mathbf 1_{|\und t(x)|=1}
	\end{equation*}
Using that the clocks are independent and \eqref{n2.6} we get
	\begin{equation*}
 E^\ga\big[ F\big]= \frac 1N\sum_{x\in D}\sum_{ y\in D}
 \Big( e^{-\delta} \delta\{
e^{-\la^*\delta A_\ga( x , y)}\la^*\delta  A_\ga( x , y)\} \Big)\le c \xi^d \delta^{2}
	\end{equation*}
We next compute the variance and using independence we get
	\begin{eqnarray*}
&&  E^\ga\Big[(F- E^\ga[F])^2\Big]= \frac 1{N^2}\sum_{x\in D, y\in D}  \delta e^{-\delta}
e^{-\la^*\delta A_\ga( x , y) }\la^*\delta A_\ga( x , y)  \le \frac 1N c \xi^d\delta^{2}
%
	\end{eqnarray*}
The Lemma  follows from the Chebishev inequality.
\qed

\medskip

\noindent
{\bf Proof of Theorem \ref{thm.n2.1.33}.}
  We fix $\bar x$  as in the hypothesis of the Theorem and we call  $\bar D=D(\bar x)$.

%
%

\medskip

Recalling \eqref{3.18.37} we write
		\begin{eqnarray}
\Big|N^{-1}\mathcal M_{D} -   \delta  v_{k,n\delta}(x)\Big|\le \big|S_N-  \delta  v_{k,n\delta}(x)|
+ \frac 1NR_{  \bar D}+ \frac 1N T_{\bar D}
	\label{3.22a}
		\end{eqnarray}
where
\begin{equation}
	\label{3.22b}
 S_N=\frac 1N \sum_{x\in \bar D}\mathbf 1_{|\und t(x)|=1},
		\end{equation}
The term $R_{  \bar D}$ has been treated in Lemma \ref{lemr.33} and the one with $T_{\bar D}$  in Lemma \ref{lemt.333}.

We  first compute $\mathbb E^\ga[S_N]$
\begin{eqnarray}
E^\ga[S_N]=e^{-\delta }\delta
\frac{|\bar D|}{N}= e^{-\delta }\delta v_{k,n\delta}(\bar x)
	\label{3.24b}
\end{eqnarray}
 Since the clocks are independent we get
	\begin{equation}
	\label{3.25}
 \mathbb E^\ga\Big[\big(S_N-\mathbb E^\ga[S_N]\big)^2\Big]=\frac 1{N^2} \sum_{x\in \bar D}\big\{P^\ga(|\und t(x)|=1)-P^\ga(|\und t(x)|=1)^2\}\le c\frac {\delta}N
		\end{equation}
Since $1-e^{-\delta }\le c \frac \delta N$, given any $\theta>0$ and using  \eqref{3.25}  we get
	\begin{eqnarray}\label{3.26}
P^\ga\Big[\big|S_N- \delta  v_{k,n\delta}(x))\big|>N^{-\theta}\Big]\le c \frac \delta N+c\frac {\delta}{N^{1-2\theta}} 	\end{eqnarray}
\qed

\section{ Proof of Theorem \ref{thma1.3bis}}
\label{sec55}

The aim in this section is  to study the time-continuum limit of the densities when first $N=(\ga^{-1}\xi)^d\to \infty$ and then $\delta\to 0$. To this end we study the process in the time interval  $\dis{[0,T]=\bigcup_{n< \delta^{-1}T}[n\delta,(n+1)\delta)}$.

\bigskip
Define
	\begin{equation}
	\label{4.2}
\mathcal X'_{N,\delta}(n) =\bigcap_{i<k}\bigcap_{x,y\in\Om_\ga}\Big\{ \big|N^{-1}\mathcal M_{ C(x) , D(y); i} (n)-   \la^* \delta N A_\ga(x,y)  v_{k,n\delta}(x) v_{i,n\delta}(y)
 \big|\le  c( N^{-\theta}+\delta^{1+\eps})\Big\}
 \end{equation}

 \begin{equation}
\label{4.4}
 \mathcal X''_{N,\delta}(n) = \bigcap_{x\in\Om_\ga}\Big\{\big|N^{-1}\mathcal M_{ D(x)}(n) -
 \delta    v_{k,n\delta}(x)
 \big|\le  c( N^{-\theta}+ \delta^{1+\eps})\Big\}, \end{equation}

\begin{lemma}
\label{lem2.1}
Let
 \begin{equation}
\label{4.8}
 \mathcal G_\ga = \bigcap_{n\delta \le T} (\mathcal X'_{N,\delta}(n)  \cap \mathcal X''_{N,\delta}(n))
 \end{equation}
 then
\begin{equation}
\label{4.9}
 P^\ga[\mathcal G_\ga] \ge 1- \delta^{-1}T c(N^{-1+2\theta} +\delta^{1+\eps})
 \end{equation}
\end{lemma}

\noindent{\bf Proof}. Observe that the sets in the curly brackets in \eqref{4.2} and \eqref{4.4} are constant on the cubes, then the lemma follows from Theorems \ref{thm.n2.1} and  \ref{thm.n2.1.33}.\qed

%
%
%
%

\medskip

We first take $N\to \infty$ and then $\delta\to 0$, thus the probability
of $\mathcal G_\ga$ is as close to 1 as we want if for any $\delta$ small enough we take $N$ sufficiently large.

\bigskip

 To underline the dependence of $v$ on $\ga$   we write
below $v^\ga_{i,n\delta}(y)$ and   we
 rewrite \eqref{n2.2} as
\begin{equation}
\label{4.5}
 v^{\ga}_{i,n\delta}(x) = \frac 1N \sum_{y\in C(x)} \mathbf 1_{\{U_{n\delta}(y)=i\}},\quad
 \und v^\ga(n\delta) = \{ v^{\ga}_{i; n\delta}(x), x\in \Om_\ga, i \in\{0,\ldots,k\}\},\quad  n\delta \in[0,T]
 \end{equation}
For $1\le i \le k$ we call
\begin{eqnarray}
\label{4.6}
&&\hskip-1cm   F_\ga(x,i;\und v^\ga(n\delta)):= \{v^{\ga}_{i-1;n\delta}(x)-v^{\ga}_{i;n\delta}(x)\}\sum_y \la^*   A_\ga(y,x) v^{\ga}_{k;n\delta}(y)
 - v^{\ga}_{k;n\delta}(x)\mathbf 1_{i=k}
 \end{eqnarray}
 \begin{eqnarray}
\label{4.7}
&&\hskip-2cm   F_\ga(x,0;\und v^\ga(n\delta)):= -v^{\ga}_{0;n\delta}\sum_y \la^*   A_\ga(y,x)  v^{\ga}_{k;n\delta}(y)
 +   v^{\ga}_{k;n\delta}(x)
 \end{eqnarray}
Recalling \eqref{n2.7.1},  \eqref{n2.7.2} and \eqref{n2.7.3}
in  $\mathcal G_\ga$ we have that for any $n\le \delta^{-1}T$ and $i\in \{0,1,..,k\}$,

\begin{eqnarray}
\label{4.10}
| v^{\ga}_{i;n\delta}(x) - \{v^{\ga}_{i;(n-1)\delta}(x)+\delta F_\ga(x,i;\und v^\ga((n-1)\delta))\}|
 \le c (N^{-\theta}+ \delta^{1+\eps})
 \end{eqnarray}
  having used that $\dis{\sum_y  A_\ga(y,x)v^{\ga}_{k;n\delta}(y)=\sum_C NA_\ga(y_C,x) v^{\ga}_{k;n\delta}(y_C)}$, with $y_C$ any point in $C$.
 \bigskip

Define $u^{\ga}(x,i;n\delta)$, $i\in \{0,1,..,k\}$, $n\delta\le T$ as the solution of the equations
     \begin{eqnarray}
\label{4.11}
&& u^{\ga}(x,i;n\delta) = u^{\ga}(x,i;(n-1)\delta) +\delta F_\ga(x,i;\und u^{\ga}((n-1)\delta))\nn \\&&
 u^{\ga}(x,i;0)= v^{\ga}_{i;0}(x)
     \end{eqnarray}
Observe that since the initial datum is constant on the cubes $C$ of the partition then also $u^{\ga}(x,i;n\delta)$ is constant on the cubes for any $n$. We thus may call
 $u_i(C,n)=u^{\ga}(x,i;n\delta)$ with any $x\in C$.

     \medskip

    \begin{lemma}

    Let  $\eps < u^{\ga}(x,i;0) <1 -\eps$ (see Definition \ref{defin3.0}), let $T>0$ and $\eps'$ such that
       \begin{equation}
\label{4.999.1}
   (1-\delta)^{\delta^{-1}T}\eps =: \eps'
       \end{equation}
    Then   $\eps' < u^{\ga}(x,i;n\delta) <1 -\eps'$ for all $x$,  $i$  and $n \le \delta^{-1}T$.
    \end{lemma}

    \noindent
    {\bf Proof.}
  Let $C$ be a cube of the basic partition, and let
    \begin{equation}
\label{4.999.1}
 K_n(C) = \sum_x \la^* A_\ga(x,y)u^{\ga}(x,k;n\delta),\quad y \in C
    \end{equation}
 observing that the right hand side does not depend on which $y$ we take in $C$.    We then have
     \begin{equation}
\label{4.999.2}
  u_i(C,n+1) - u_i(C,n) = K_n(C)\Big( u_{i-1}(C,n) -u_i(C,n)\Big)\delta - \delta u_i(C,n) \mathbf 1_{i=k}
    \end{equation}
with $u_{-1}(C,n)\equiv 0$.
We call
        \begin{equation}
\label{4.999.4}
 \theta(n) =\min_i u_i(C,n),\quad  \Theta(n) =\max_i u_i(C,n)
    \end{equation}

     \centerline{\em The upper bound}

     \medskip
\noindent
Let $i$ such that $\Theta(n+1)=u_i(C,n+1)$ and $j$ such that $\Theta(n)=u_j(C,n)$.  By \eqref{4.999.2}
bounding $u_{i-1}(C,n)\le u_j(C,n)$ and dropping the last term $-\delta u_{i}(C,n)\mathbf 1_{i=k}\le 0$ we get
           \begin{eqnarray}
\label{4.999.8}
  u_i(C,n+1)- u_j(C,n) &\le&    u_i(C,n)- u_j(C,n) + K_n(C)[u_{j}(C,n)-u_{i}(C,n)]\,\delta
  \nn\\&\le& - (u_j(C,n)- u_i(C,n))\Big(1-K_n(C)\,\delta\Big)
           \end{eqnarray}
It follows by induction from \eqref{4.999.8}  that   $K_n(C)\le \la^*$ for all $n$. Thus the right hand side of \eqref{4.999.8} is negative
for $\delta$ small enough and therefore $\Theta (n) \le \Theta (0) \le \eps$ which proves the upper bound.
 \medskip

 \centerline{\em The lower bound}

\noindent
Let $i$ such that $\theta(n+1)=u_i(C,n+1)$ and $j$ such that $\theta(n)=u_j(C,n)$.  We prove below that for $\delta$ small enough
       \begin{equation}
\label{4.999.5}
  u_i(C,n+1)\ge u_j(C,n)
   - u_j(C,n)\delta\quad \text{ namely }\qquad \theta(n+1)\ge \theta(n)(1-\delta)
    \end{equation}
From \eqref{4.999.5} we get
for any $n$ such that $n\delta^{-1} \le T$
       \begin{equation}
\label{4.999.6}
 \theta(n)\ge (1-\delta)^{\delta^{-1}T}\eps =: \eps'
    \end{equation}
having used Definition \ref{defin3.0}.

    \medskip

\noindent
{\em Proof of  \eqref{4.999.5}}.
Recalling the definitions of $i$ and $j$ from \eqref{4.999.2} we have
           \begin{eqnarray}
\label{4.999.7bis}
  u_i(C,n+1)- u_j(C,n) &=&    u_i(C,n)- u_j(C,n) + K_n(C)[u_{i-1}(C,n)-u_{i}(C,n)]\,\delta
  \nn\\&-&\delta u_{i}(C,n) \mathbf 1_{i=k}
           \end{eqnarray}
 We bound from below   $ u_{i-1}(C,n)  \ge u_j(C,n)$.  We also write the last term as   $[(u_{i}(C,n)-u_j(C,n)) +u_j(C,N)]\,\delta$ having bounded $ \mathbf 1_{i=k}\le 1$ and get
            \begin{equation}
\label{4.999.7}
\theta(n+1) - \theta(n) \ge \{ [u_i(C,n)- u_j(C,n)] \Big(1 - K_n(C)\delta -\delta\Big)\}
-u_j(C,n)]\,\delta
%
     \end{equation}
For  $\delta$ small enough the curly bracket term is positive and \eqref{4.999.5} is proved. \qed
    \medskip

 Let
 \begin{equation}
\label{4.12}
 \|\und v^\ga(n\delta)-\und u^\ga(n\delta)\| := \max_{x,i}|v^\ga_{i;n\delta}(x)-u^{\ga}(x,i;n\delta) |
 \end{equation}

 \begin{lemma}
 \label{lem4}
 In $\mathcal G_\ga$ we have
 \begin{equation}
\label{a4.17a}
\sup_{n:n\delta<T}\|\und v^\ga(n\delta)-\und u^\ga(n\delta)\| \le   c T\delta^{-1} (N^{-\theta}+ \delta^{1+\eps})
\end{equation}
which vanishes in the limit $N\to \infty$ and then $\delta \to 0$.
\end{lemma}
\medskip

\noindent{\bf Proof.}
There is $c$ so that
\begin{equation}
\label{4.13}
\max_{x,i}| F_\ga(x,i;\und v^\ga(n\delta)-F_\ga(x,i;\und u^\ga(n\delta)| \le c \|\und v^\ga(n\delta)-\und u^\ga(n\delta)\|
\end{equation}
and therefore by \eqref{4.10}
\begin{equation}
\label{4.14}
\|\und v^\ga(n\delta)-\und u^\ga(n\delta)\| \le c (1+\delta)\|\und v^\ga((n-1)\delta)-\und u^\ga((n-1)\delta)\| + c(N^{-\theta}+ \delta^{1+\eps})
\end{equation}
By iteration
\begin{equation}
\label{4.15}
\|\und v^\delta(n)-\und u^\delta (n)\| \le \sum_{m=2}^n [c (1+\delta)]^m  (N^{-\theta}+ \delta^{1+\eps})
\end{equation}
Since
\begin{equation}
\label{4.16}
 [c (1+\delta)]^{T\delta^{-1}} \le  c'\equiv c'(T)
\end{equation}
and $n\le T\delta^{-1}$, we get
\begin{equation}
\label{4.17}
\|\und v^\ga(n\delta)-\und u^\ga(n\delta)\| \le   c' T\delta^{-1} (N^{-2\theta}+ \delta^{1+\eps})
\end{equation}
which concludes the prove the Lemma.\qed

\bigskip
We next define $w^{\ga}(x,i;t)$, $t\in[0,T]$, which is the time continuous analogue of $u^\ga$, as the solution of
\begin{equation}
\label{4.18}
\frac{d}{dt}w^{\ga}(x,i;t)= F_\ga(x,i;\und w^{\ga}(t)),\quad
\und w^{\ga}(t) = \big\{w^{\ga}(x,i;t), x\in\Om_\ga, i\in\{0,..,k\}\big\}
\end{equation}
with initial condition $\und w^{\ga}(0)=\und v^\ga(0)$.  Then
\begin{eqnarray}
\label{4.19}
&&| w^{\ga}(x,i;n\delta) - \{w^{\ga}(x,i;(n-1)\delta) +\delta F_\ga(x,i;\und w^{\ga}((n-1)\delta)\}| \le c\delta^2
 \end{eqnarray}
and
\begin{equation}
\label{4.20}
\|\und  w^{\ga}(n\delta)-\und u^\ga(n\delta)\| \le   c' T\delta^{-1} \delta^2
\end{equation}
Therefore
\begin{equation}
\label{4.21}
\|\und  w^{\ga}(n\delta)-\und v^\ga(n\delta)\| \le   c T\Big(\delta +\delta^{-1}N^{-\theta}+ \delta^\eps\Big),\qquad \forall n\le T\delta^{-1}
\end{equation}

\bigskip
\noindent{\bf Proof of Theorem \ref{thma1.3bis}}.
We call
\begin{equation}
\label{4.22}
 \om^\ga(r,i;t) := w^{\ga}(\ga^{-1} r,i;t),\qquad r\in \Om^*_\ga
\end{equation}
and we observe that $ \om^\ga$ satisfies
\begin{equation}
\label{4.22bis}
\frac{d}{dt}\om^{\ga}(r,i;t)= f_\ga(r,i;\und \om^{\ga}(t)),\quad
\und \om^{\ga}(t) = \{\om^{\ga}(r,i;t), r\in\Om^*_\ga, i=0,..,k\}, t\in[0,T]
\end{equation}
with initial condition $\om^{\ga}(r,i;0)= w^{\ga}(r,i;0)$ and 
\begin{eqnarray}
\label{4.23}
&&\hskip-2cm   f_\ga(r,i;\und \om^{\ga}(t)):= \{
\om^{\ga}(r,i-1;t)
-\om^{\ga}(r,i;t)\} \sum_{C^\xi}  \la^*  N A_\ga(\ga^{-1}r',\ga^{-1}r)\om^{\ga}(r',k;t)
 \end{eqnarray}
 where $r'$ is any point in $C^\xi$. We have used that $w^{\ga}$ is constant on the cubes $C^{\ga,\xi}$ as well as $A_\ga$. We now observe that (recall \eqref{nuovo})
 	\begin{equation}
	\label{4.34}
\lim_{\ga \to 0} \ga^{-d}A_\ga(\ga^{-1}r',\ga^{-1}r)=\xi^d A_\xi(r',r)
	\end{equation}
 Thus, since $N=(\ga^{-1}\xi)^d$ we get that
\begin{equation}
	\label{4.35}\lim_{\ga\to 0} \|\und\om^{\ga}(t)-\und \varphi_\xi(t)\|=0, \quad
	\und\varphi_\xi(t) = \{\varphi_\xi(r,i,t), r\in\Omega, i\in\{0,\ldots,k\}\}	
	\end{equation}
 From Lemma \ref{lem2.1},  \eqref{4.21} and \eqref{4.35} we then get

 	\begin {equation*}
\lim_{\delta\to 0}\lim_{\ga\to 0}	P^\ga\Big[\|\und \varphi_\xi(t)-v^\ga(n\delta)\| \le   T\big(\delta +\delta^{-1}N^{-\theta}+ \delta^\eps\big)\Big]
=1
	\end{equation*}
which proves the Theorem. \qed

\bigskip

\section{Stronger version of Theorems \ref{thma1.3bis} and \ref{thm2.3a}}
\label{sub5.1}

In this Section we prove Theorem \ref{thm2.3a} together with a stronger version of Theorem \ref{thma1.3bis}.

We introduce some new notation and definitions
 besides  those in the previous sections.

 We fix a time $t=n_0\delta_0$ for some $n_0$ and $\delta_0$. Since we consider the parameter $\delta$ of the form $2^{-n_2}$ with $n_2\in\mathbb N$ we then have that for any $\delta<\delta_0$ there is $m$ so that $t=m\delta$.

 We define $v_{C,i;t}$ with $C\in\pi_\xi$, $\xi=2^{-n_3}$, $n_3\in\mathbb N$, $i\in\{0,1,..,k\}$.
 \begin{equation}
	\label{5.1}
v_{C,i;t}= \frac 1{|C|}\sum_{r\in C}\mathbf 1_{U_t(r)=i},\qquad |C|=\xi^d				\end{equation}
where by an abuse of notation we have called $U(r)$, $r=\ga x$, the potential $U(x)$;
$U_t(r)$ is the potential at time $t$.

So far we have studied the one-body correlation functions.
In the next theorem we study  the many body space-time correlations, namely  the law $\mathcal P^{\ga,\delta,\xi}_{\und v} $ of the
finite dimensional distribution $\und v= \big(v_{C_\ell,i_\ell;t_\ell}, \ell=1,2,...,m\big)$  in the limit first $\ga\to 0$ then $\delta\to 0$ and finally $\xi \to 0$.

\medskip

\begin{thm}
\label{teo5.1}
With the above notation
			\begin{equation}
	\label{5.3}
\lim_{\xi\to 0}\,\,\lim_{\delta\to 0}\,\,\lim_{\ga\to 0} P^{\ga,\delta,\xi}_{\und v} =P_{\und v}
			\end{equation}
where $P_{\und v}$ is supported by $\und w=(w_1,..,w_m)$ with
	\begin{equation}
	\label{5.4}
w_\ell=\frac 1{|C_\ell|}\int_{C_\ell} dr \rho(r,i_\ell;t_\ell),\qquad \ell=1,2,...,m
			\end{equation}
 $\rho$ is the solution of \eqref{1.2a}-\eqref{1.2c}.
\end{thm}

\noindent{\bf Proof.} Recalling the definition of $\om^\ga$ in \eqref{4.22} we write
\begin{equation}
	\label{5.5}
v_{C_\ell,i_\ell;t_\ell}=\{v_{C_\ell,i_\ell;t_\ell}-\om^\ga(r,i_\ell;t_\ell)\}+\om^\ga(r,i_\ell;t_\ell)\qquad r\in C_\ell
			\end{equation}
In the set
$ \mathcal G_\ga $ defined in \eqref{4.8} the curly bracket is uniformly bounded by $c( N^{-\theta}+\delta^{\eps})$  and by Lemma \ref{lem2.1} $P^\ga(\mathcal G_\ga)\to 1$ as $\ga\to 0$. Thus in the limits $N\to \infty$ and then $\delta\to 0$  the curly bracket  goes to 0.

From \eqref{4.35} $\om^\ga(r,i_\ell;t_\ell)$ converges as $\ga\to 0$ to $\varphi(r,i_\ell;t_\ell)\equiv \varphi_\xi(r,i_\ell;t_\ell)$ solution of
\eqref{n2.3.3.0},\eqref{n2.3.3.0bis}, \eqref{n2.3.3.0ter} where the convolution term $A\equiv A_\xi$ is defined in \eqref{nuovo} and
 \begin{equation}
	\label{5.6}
\lim_{\xi\to 0}A_\xi(r,r')= \lim_{\xi\to 0} \,\, \frac{1}{\xi^{2d}}\int_{  C^{\xi}(r)}dr_1\int_{ D^{\xi}(r')}dr'_1
 J(r_1,r_1')	=J(r,r')\end{equation}

 From \eqref{5.6} it is easy to prove the
 $$\lim_{\xi\to 0} \varphi_\xi(r,i_\ell;t_\ell)=w(r,i_\ell;t_\ell)$$
and the Theorem is proved.\qed

\medskip

%
%

{\bf{Positive, real valued times.}}
Even though the set $\mathcal T$ is dense in $[0,T]$, yet it sounds non physical to restrict times to $\mathcal T$.  The problem can be fixed easily using a variable time mesh.

To explain the idea we refer first to the simpler case of a single time $t\in [0,T]$ as in Theorem \ref{thma1.3bis}.  Suppose $t\notin \mathcal T$. We then consider a mesh \ $\delta \in \{2^{-n}t, n \in \mathbb N\}$, and similarly a second mesh $\delta' \in\{ 2^{-n}(T-t), n \in \mathbb N\}$.  We can then use the proof of Theorem  \ref{thma1.3bis} in $[0,t]$ where the mesh is $\delta$ and again the proof of Theorem  \ref{thma1.3bis} in $[t,T]$ with the mesh $\delta'$.  The extension to the case of Theorem \ref {thm2.3a} is similar. We have $m$ times $0<t_1,\cdots < t_m< T$ we then consider a mesh $\delta_1 \in \{2^{-n}t_1, n \in \mathbb N\}$, $\dots$,  $\delta_{m+1} \in \{2^{-n}t_{m}, n \in \mathbb N\}$,
and use
the proof of Theorem  \ref{thma1.3bis} in each one of the above time intervals.

\section{Extensions}
\label{sec5}

In this Section we   study
the macroscopic limit of other infection/recovery models.

\medskip

\subsection{Additional recovery jumps}
\label{sub5.3}
\medskip

Here we consider the case where also individual with potential $i<k$ may recover i.e. the individual at
 site $x$ with $U(x)=i$ recovers at rate $\la_i$ and the potential after recovery becomes $U(x)= 0$.

The macroscopic equations are:
 \begin{eqnarray}
	\label{5.7a}
\frac{d}{dt}\rho_i(r,t)=-\la_i\rho_i(r,t)+\{ \rho_{i-1}(r,t) -\rho_i(r,t)\} \int_{\Om}dr'\la^*\rho_k(r',t) J(r',r)
	\end{eqnarray}
when  $1\le i \le k-1$; for $i=k$
\begin{equation}
	\label{5.7aa}
\frac{d}{dt}\rho_k(r,t)= - \rho_k(r,t)+ \rho_{k-1}(r,t)\int_{\Om}dr'\la^*\rho_k(r',t)J(r',r)
	\end{equation}
and  for $i=0$:
\begin{equation}
	\label{5.8aa}
\frac{d}{dt}\rho_0(r,t)= \sum \lambda_i \rho_i (r,t) +  \rho_k(r,t)- \rho_0(r,t)\int_{\Om}dr'\la^*\rho_k(r',t) J(r',r)
	\end{equation}

The proof is similar to the proof of Theorem \ref{teo5.1} and omitted.

\medskip

\subsection{The excitatory-inhibitory
network model}
\label{sub5.4}

\medskip
Referring to a neural network here we consider excitatory and inhibitory neurons; both neurons have a potential in $\{0,..,k\}$. When a excitatory neuron with potential $k$ fires the  potentials of all the other neurons with potential $<k$  increase by 1. Similarly when an inhibitory neuron with potential $k$ fires the  potentials of all the other neurons with potential $>0$ and $<k$ decrease by 1. The rates of firing are $\la^*_1 k J^{(j)}_\ga(x,y)$, $j=1$ for excitatory and $j=2$ for  inhibitory. Besides that neurons with potential $k$ decay at rate 1 to a state with potential 0.

For this model we derive the following macroscopic equations. We denote by $\rho_1(r,i;t)$ the limit macroscopic density at position $r$ and time $t$ of the excitatory neurons with potential $i$ and by $\rho_2(r,i;t)$ the limit density of the inhibitory neurons with potential $i$.

 \begin{eqnarray}
	\nn
\frac{d}{dt}\rho_j(r,i;t)&&=\{ \rho_j(r,i-1;t) -\rho_j(r,i;t)\} \int_{\Om}dr'\la_1^*\rho_1(r',k;t) J^{(1)}(r',r)
\\&&+ \{ \rho_j(r,i+1;t) -\rho_j(r,i;t)\} \int_{\Om}dr'\la_2^* \rho_2(r',k;t) J^{(2)}(r',r)
	\end{eqnarray}
when  $1\le i \le k-1$; for $i=k$
\begin{equation}
	\label{6.14a}
\frac{d}{dt}\rho_j(r,k;t)= - \rho_j(r,k;t)+ \rho_j(r,k-1;t)\int_{\Om}dr'\la_1^*\rho_1(r',k;t)J^{(1)}(r',r)
	\end{equation}
for $i=0$
\begin{eqnarray}
	\nn
\frac{d}{dt}\rho_j(r,0;t)&&=  \rho_j(r,k;t)+ \rho_j(r,1;t)\int_{\Om}dr'\la_2^*\rho_2(r',k;t)J^{(2)}(r',r)
\\&&- \rho_j(r,0;t)\int_{\Om}dr'\la_1^*\rho_1(r',k;t)J^{(1)}(r',r)
\label{6.14aa}
\end{eqnarray}
Convergence is in the sense of the finite dimensional distributions as in Theorem \ref{teo5.1}.

%
\medskip

\subsection{\bf General microscopic model.}
Place at each site $x\in\mathbb{Z}^d$ a finite-state, continuous-time Markov chain $U(x,t)$ with state space $S=\{0,\ldots,k\}$.
 For any pair of states $i$ and $j$ there is an intrinsic transition rate from $i$ to $j$, denoted as $g^{ij}_\gamma(x)$, dependent on the scaling parameter $\gamma$ and location $x$.  Any other site $y$ at any state $U(y,t)=l\in\{0,\ldots,k\}$ will have an additive effect $\lambda_{ijl}J^{ijl}_\gamma (y,x)$ on the transition rates at $x$ from $i$ to $j$.  Then together, for any $i\neq j\in S$, the time-dependent transition rate from $i$ to $j$, denoted $q_{ij}(x,t)$, is given by
\[
q_{ij}(x,t) = g^{ij}_\gamma (x) + \sum_l \lambda_{ijl}\sum_{y\in\mathbb{Z}^d}J^{ijl}_\gamma (y,x) \delta(U(y,t)-l),
\]
where $\lambda_{ijl}\geq 0$, $g^{ij}_\gamma (x)=g^{ij}(\gamma x)$, and $J^{ijl}_\gamma (x,y)=\gamma^d J^{ijl}(\gamma x, \gamma y)$.  The functions $g^{ij}(x)$ and $J^{ijl}(x,y)$ describe the intrinsic transition rates and the site-to-site interactions in the scaling limit.  We take the assumptions on $J^{ijl}$ to be the same as before, and we assume $g^{ij}$ to be continuous and bounded over $x$.

\medskip\noindent{\bf The scaling limit.} In the scaling limit $\gamma \rightarrow 0$, the state of the system is described by local state distributions $\vec{v}(x,t) = (v_1(x,t),\ldots, v_N(x,t))$ which vary continuously in $x\in\mathbb{R}^d$, and evolve in time according to macroscopic equations, for $0\leq i \leq k$,
\begin{align}
v_i' = \sum_{j\neq i} \left[g^{ji} + \sum_l \lambda_{jil}R_{jil}(x) \right]v_j-\left[g^{ij} +\sum_l \lambda_{ijl}R_{ijl}(x)\right]v_i
\end{align}
with
\begin{align}
R_{ijl}(x) = \int J^{ijl}(y,x)v_k(y)dy.
\end{align}

\medskip\noindent{\bf Obtaining the original, generalized contact process.} The original, generalized contact process can be recovered by choosing particular $g^{ij}$ and $J^{ijl}$.  Let $g^{k0}=1$ and all other $g^{ij}=0$ for $i\neq j$.  Next, for $0\leq i < k$, let $\lambda^{i,i+1,k}=\lambda^*$ and all other $\lambda^{ijl}=0$.  Let $J^{i,i+1,k}$ be defined as for the generalized contact process.  Then the microscopic and macroscopic equations are the same as for the generalized contact process.

\medskip\noindent{\bf Adding new features.} The E-I neuron model from subsection \ref{sub5.4} can be recovered by letting $S=\{0,\ldots,k\} \times \{0,\ldots,k\}$, where the first coordinate corresponds to the E neuron voltage, and second coordinate the I neuron voltage at a site.  Excitatory interactions will correspond to setting the terms $\lambda^{(i,j),(i+1,j),(k,l)}$, $\lambda^{(i,j),(i,j+1),(k,l)}$, as well at the associated $J$ terms, and inhibitory interactions will correspond to setting the terms $\lambda^{(i,j),(i-1,j),(l,k)}$, $\lambda^{(i,j),(i,j-1),(l,k)}$, as well as their associated $J$ terms.

An external drive can be added by making nonzero the terms $g^{(i,j),(i+1,j)}$ and $g^{(i,j),(i,j+1)}$.  A neuronal leak could be modeled at least approximately by setting $g^{(i,j),(i-1,j)}$ and $g^{(i,j),(i,j-1)}$.

 The proofs of all these extensions follow the same steps as in the proof of Theorem \ref{teo5.1} with  new clock processes associated to $g^{ij}_\gamma (x) $ and
$\lambda_{ijk}J^{ijk}_\gamma (y,x)$. As before all these clocks are independent and we can repeat the arguments used in Sections 3-7.

\appendix

\section{Clusters and branching processes}
\label{app.A}

In this appendix we will prove Theorem
\ref{thm?}.  For completeness we first recall some definitions and notation.

A configuration $\{m(z,z'),m(z''\}$ is the set of arrows and marked points with their multiplicity which are described respectively by integer valued functions $m(z,z')$ and $m(z'')$, $m(z,z')=0$ if the arrow $(z,z')$ is absent and $m(z'')=0$ if $z''$ is not a marked point.  Then the probability of a configuration is
\begin{equation}
\label{nA1}
P(\{m(z,z'),m(z''\}) = \{\prod_{(z,z')} e^{-\la^*\delta A( z,z')}
\frac{[\la^* \delta A( z, z')]^{m(z,z')}}{m(z,z')!}\}
\{ \prod_{z''} e^{-\delta}
\frac{\delta^{m(z'')}}{m(z'')!}\}
\end{equation}
where $P=P^{A_\ga,\ga}$ and $A( z,z')= A_\ga(z,z')$.

Let $V=V(\mathbf C)$ be the set of  arrows in $\mathbf C$, then
 a cluster $\mathbf C$ is a maximal connected set of arrows  plus the specification of the multiplicity $m(z,z')$ of the arrows
and the multiplicity $m(z'')$ of the points $z''\in T_V$, where $T_V$ is the union of starting and endpoints of the arrows in $V$.  Maximality means that there is no arrow starting from $T_V^c$
and ending at $T_V$.  Then
\begin{equation}
\label{nA2}
P(\mathbf C) = \{\prod_{(z,z')\in V} e^{-\la^*\delta A( z,z')}
\frac{[\la^* \delta A( z, z')]^{m(z,z')}}{m(z,z')!}\}
\{ \prod_{z''\in T_V} e^{-\delta}
\frac{\delta^{m(z'')}}{m(z'')!}\}\{\prod_{z\in T_V^c,z'\in T_V} e^{-\la^*\delta A( z,z')}\}
\end{equation}

With reference to \eqref{n2.6.1.1} 
we fix an arrow $(x,y)$  and write
	\begin{equation}
\label{nA3aa}
\sum_{j\ge 1}\delta^{-aj} \sum_{\mathbf C_j \ni (x,y) }P\big[\mathbf C_j\big]=S
\end{equation}
where
\begin{equation}
\label{nA3}
S:=\sum_{\mathbf C \ni (x,y)}P \big[\mathbf C\big]\delta^{-a|\mathbf C|},\quad a \in (0,1)
\end{equation}
where $|\mathbf C|$ is the number of elements in $\mathbf C$, namely 
\begin{equation}
\label{nA4}
|\mathbf C|=\sum_{(z,z)\in V} m(z,z') +\sum_{z'' \in T_V}m(z'')
\end{equation}
The purpose is to prove 
\begin{equation}
\label{nA8}
S \le c \delta^{1-a-2\eps}\frac 1N
\end{equation}
which then implies \eqref{n2.6.1.1}.

We perform the sum over $\mathbf C \ni (x,y)$ by first summing over the multiplicities and get
\begin{eqnarray}
\label{nA5}
S &=& \sum_{V\ni (x,y)}\{\prod_{(z,z')\in V} e^{-\la^*\delta A( z,z')}[e^{\la^*\delta^{1-a} A( z,z')}-1]\}\{ \prod_{z''\in T_V} e^{-\delta}e^{\delta^{1-a}}
\}\nn\\&\times&\{\prod_{z\in T_V^c,z'\in T_V} e^{-\la^*\delta A( z,z')}\}
 \end{eqnarray}
Since $A(z,z') \le c \frac{\xi^d}N$,
 \begin{equation}
\label{nA6}
|e^{\la^*\delta^{1-a} A( z,z')}-1| \le c\la^*\delta^{1-a}\frac{\xi^d}{N},\;
e^{-\delta}e^{\delta^{1-a}}\le 1+ c\delta^{1-a},\;e^{-\la^*\delta A( z,z')} \le 1
\end{equation}
 We also split $\delta^{1-a}= \delta^{1-a-2\eps}\delta^{2\eps}$,
 $\eps>0$ such that $ 1-a -2\eps>0$, and get
\begin{eqnarray}\nn
S &&\le \delta^{1-a-2\eps}
\sum_{V\ni (x,y)} \{\prod_{(z,z')\in V}\Big( c\la^*\delta^{ 2\eps} \frac{\xi^d}N\Big)\}
\{ \prod_{z''\in T_V} (1+ c\delta^{1-a})\}
\\&&\le \delta^{1-a-2\eps}
\sum_{V\ni (x,y)} \{\prod_{(z,z')\in V}\Big( c'\la^*\delta^{ 2\eps} \frac{\xi^d}N\Big)\}
\label{nA7}
\end{eqnarray}
with $c' > c(1+c\delta^{1-a})$.
The first factor because there is at least one arrow (namely $(x,y)$ which starts from $x$). We will prove that the term multiplying $\delta^{1-a-2\eps}$ is bounded by $c/N$ and  thus prove \eqref{nA8}.
The proof will exploit the  branching structure of $V$.

\medskip

We call $x$ the root of the branching,
$(x,z^1_{1}),..,(x,z^1_{n_1})$, $z^1_1=y$, the arrows which start from the root $x$ ($n_1\ge 1$ because $V\ni(x,y)$)
and $z^1_{1} ,.., z^1_{n_1} $  the nodes of the first generation. From each node
$z^1_{i}$ of the first generation may or may not start new arrows: if no arrow starts from all the nodes $z^1_i$ then the branching ends, otherwise we call $z^2_{1} ,.., z^2_{n_2} $ the nodes which are the endpoints of the new arrows: these are the nodes of the second generation.
Notice that there may be arrow which go back to $x$, in that case that arrow will not produce descendants  because they are already included in the arrows of the first generation.
Analogously we call $\{z^{i}_1, ..,z^{i}_{n_{i}}\}$ the endpoints of the arrows starting from nodes of the $(i-1)$-th generation. The branching ends when no arrow starts from the nodes of the last generation. 

In terms of the branching the configuration are described by the following parameters:
\begin{equation*}
k;\quad (n_i, 1\le i\le k);\quad (z^{i}_j,1\le i\le k, 1\le j\le n_i);\quad (R^{i}_j,1\le i\le k-1, 1\le j\le n_i)
\end{equation*}
which must satisfy
	\begin{equation}
	\label{app}
\sum_{j=1}^{n_i} R^{i}_j= n_{i+1}
\end{equation}
where

$\bullet$ $k\ge 1$ is the number of generations,

$\bullet$  $z^{i}_j$ are the position of the nodes,

$\bullet$ $R^{i}_j$ is the number of arrows which start from $z^{i}_j$.

\medskip
Writing \eqref{nA7} in terms of the branching we get
\begin{equation*}
S \le \delta^{1-a-2\eps}
\sum_{k\ge 1} \sum_{(n_i)}\sum_{(z^{i}_j)}\sum_{(R^{i}_j)}\Big( c'\la^*\delta^{2 \eps} \frac{\xi^d}N\Big)^{n_1}
\{\prod_{i=1}^{k-1}\Big( c'\la^*\delta^{2 \eps} \frac{\xi^d}N\Big)^{R^i_1+\dots+R^i_{n_i}}\}
\end{equation*}
By  \eqref{app} and since $z^{i}_j$ is the endpoint of an arrow, so that it may have at most $ N$ values, then the sum over $(z^{i}_j)$  is bounded by $N^{n_1-1+n_2+..+n_k}$ (recall that $z^1_1=y$) we get
\begin{equation*}
\sum_{(z^{i}_j)}\Big(  \frac{1}N\Big)^{n_1}
\{\prod_{i=1}^{k-1}\Big(  \frac{1}N\Big)^{R^i_1+\dots+R^i_{n_i}}\}\le \frac 1N
\end{equation*}
We thus have
\begin{eqnarray*}
S&& \le \frac 1N\delta^{1-a-2\eps}
\sum_{k\ge 1} \sum_{(n_i)}\sum_{(R^{i}_j)}\Big( c'\la^*\delta^{2 \eps} \xi^d\Big)^{n_1}
\prod_{i=1}^{k-1}\Big( c'\la^*\delta^{2 \eps} \xi^d\Big)^{R^i_1+\dots+R^i_{n_i}}
\\&&\le \frac 1N\delta^{1-a-2\eps}
\sum_{k\ge 1} \sum_{(n_i)}[c'\la^*\xi^d]^{n_1+..n_k}\sum_{(R^{i}_j)} \delta^{2 \eps n_1} 
\prod_{i=1}^{k-1}\Big( \delta^{2 \eps} \Big)^{R^i_1+\dots+R^i_{n_i}}
	\end{eqnarray*}
We first estimate the sum on $(R^{k-1}_j)$ which satisfy \eqref{app} and get
	\begin{equation}
	\sum_{(R^{k-1}_j)} 
	(\delta^{2 \eps})^{R^{k-1}_1+\dots+R^{k-1}_{n_{k-1}}}
	\le \Big[\sum_{m=0}^{\infty} (\delta^{2 \eps})^m\Big]^{n_{k-1}}\le (1+c\delta^{2\eps})^{n_{k-1}}
	\label{app1}
	\end{equation}
writing $\delta^{2\eps}=\delta^\eps\delta^\eps$ and using \eqref{app} we get
\begin{equation}
	\sum_{(R^{k-2}_j)} 
	(\delta^{2 \eps})^{R^{k-2}_1+\dots+R^{k-2}_{n_{k-2}}}
	\le \delta^{\eps n_{k-1}}\Big[\sum_{m=0}^{\infty} ( \delta^{ \eps} )^m\Big]^{n_{k-2}}\le  \delta^{\eps n_{k-1}}(1+c\delta^{\eps})^{n_{k-2}}
	\label{app2}
	\end{equation}
The other sums on $(R^{i}_j)$ with $i<k-2$ are estimated as in \eqref{app2} getting 
\begin{eqnarray}\nn
S &&\le \frac 1N\delta^{1-a-2\eps}
\sum_{k\ge 1} \sum_{(n_i)}[c'\la^*\xi^d]^{n_1+..n_k}[\delta^\eps(1+c\delta^\eps)]^{n_1+..+n_{k-1}}
\\&&= \frac 1N\delta^{1-a-2\eps}\sum_{k\ge 1} \sum_{(n_i)}[c'\la^*\xi^d\delta^{\eps/2}]^{n_1+..n_k}[\delta^{\eps/2}(1+c\delta^\eps)]^{n_1+..+n_{k-1}}
\nn
\\&&\le
\frac 1N\delta^{1-a-2\eps}\sum_{k\ge 1}\big[ \sum_{n=1}^\infty [\delta^{\eps/2}(1+c\delta^\eps)]^{n}\big]^k\le \frac 1N\delta^{1-a-2\eps}\sum_{k\ge 1}(c\delta^{\eps/2})^k
\label{app3}\end{eqnarray}
where for $\delta$ small enough we have bounded $c'\la^*\xi^d\delta^{\eps/2}<1$.

Thus \eqref{nA8} follows from \eqref{app3}\.\qed

\medskip

The proof of \eqref{n2.6.1.2} is reduced to that of \eqref{n2.6.1.1} when we write
	\begin{equation}
	\label{app5}
	\text{l.h.s. of \eqref{n2.6.1.2}} =\frac 1N \sum_x\sum_{m(x)\ge 1} e^{-\delta}\frac {\delta^{m(x)}}{m(x)!}+ \frac 1N\sum_{x,y}S_{x,y}
 	\end{equation}
where $S_{x,y}\equiv S$ as in \eqref{nA3}, having made explicit its dependence  on $x,y$.

The first term in \eqref{app5} is the contribution of clusters with only the marked point $x$, the other clusters give rise to the second term in \eqref{app5}. \eqref{n2.6.1.2}then follows from \eqref{app5} and \eqref{nA8}.\qed

  \bigskip

\medskip\noindent{\bf Acknowledgements.} The work of JLL was supported in part by the AFOSR.  The work of LC was supported by a grant from
the Simons Foundation (691552, LC).

\bibliographystyle{amsalpha}

\begin{thebibliography}{99}

\bibitem{CL} L. Chariker, J.L. Lebowitz (2021) {\it {Time evolution and stationary states of a generalized mean-field contact process}} Journal of Statistical Mechanics Theory and Experiment, 023502.

\bibitem{DGLP} A. De Masi, A. Galves, E. Lockerbach, E. Presutti  (2015) {\it{ Hydrodynamic
Limit for Interacting Neurons}}. Journal of Statistical Physics {\bf 158} 866-902.

\bibitem{DP} A. De Masi, E. Presutti (1991) {\it {Mathematical models for hydrodynamic limits}}. Lecture Notes in Mathematics, 1501.

\bibitem{DOPT} A. DeMasi, E. Orlandi, E. Presutti, and L. Triolo (1994) {\it{ Glauber evolution with Kac potentials. I.Mesoscopic and macroscopic limits, interface dynamics}}. Nonlinearity {\bf 7} 1-67.

\bibitem{DOR}A. Duarte, G. Ost, A.A. Rodriguez (2015). {\it {Hydrodynamic Limit for Spatially
Structured Interacting Neurons}}. Journal of Statistical Physics {\bf 161} 1163–1202.

\bibitem{Gerstner}
Wulfram Gerstner and Werner Kistler.
\textit{Spiking Neuron Models: Single Neurons, Populations, Plasticity}.
Cambridge University Press, Cambridge, UK, 2002.

\bibitem{GL} G. Giacomin and J.L. Lebowitz (1997). {\it { Phase segregation dynamics in particle systems with long-range interactions I: macroscopic limits}}. J. Stat. Phys. {\bf 87} 37–61.

\bibitem{GLP} G. Giacomin, J. L. Lebowitz, E. Presutti (1999). {\it {
Deterministic and stochastic hydrodynamic equations arising from simple microscopic model systems
Stochastic partial differential equations: six perspectives}}. Math. Surveys Monogr. AMS (RI), {\bf 64}, 107-152.

\bibitem{Goychuk}
I. Goychuk, A. Goychuk (2015) {\it {Stochastic Wilson-Cowan models of neuronal network dynamics with memory and delay}}.
New J. Phys. {\bf 17} 045029

\bibitem{Kandel}
E. R. Kandel, J. H. Schwartz (1985) {\it {Principles of Neural Science}}, New York:Elsevier.

\bibitem{KL} C. Kipnis, C. Landim {\it {Scaling limits of interacting particle systems}}, Springer 1991, ISSN 0072-7830

\bibitem{Liggett}
 T.M. Liggett.
 \textit{Stochastic Interacting Systems: Contact, Voter and Exclusion Processes}. Springer-Verlag Berlin Heidelberg, 1999.

\bibitem{Marro}
J. Marro, R. Dickman.
\textit{Nonequilibrium Phase Transitions in Lattice Models}
Cambridge University Press, Cambridge, 1999.

\bibitem{Spohn} H. Spohn {\it {Large scale dynamics of interacting particles}}. Texts and Monographs in Physics. Springer-Verlag, Berlin, 1991.

\bibitem{Zankoc}
C. Zankoc, T. Biancalani, D. Fanelli, and R. Livi.
\textit{Diffusion approximation of the stochastic Wilson-Cowan model.}
Chaos, Solitons and Fractals, {\bf 103}:504-512, 2017.

\end{thebibliography}

 \end{document}